\documentclass{amsart}
\usepackage{hyperref}

\newtheorem{theorem}{Theorem}[section]
\newtheorem{lemma}[theorem]{Lemma}
\newtheorem{corollary}[theorem]{Corollary}

\newtheorem{hypo}[theorem]{Hypothesis {\bf H.}\hspace*{-0.6ex}}

\newcommand{\R}{{\mathbb R}}
\newcommand{\N}{{\mathbb N}}
\newcommand{\Z}{{\mathbb Z}}
\newcommand{\C}{{\mathbb C}}
\newcommand{\M}{{\mathbb M}}

\newcommand{\nn}{\nonumber}
\newcommand{\be}{\begin{equation}}
\newcommand{\ee}{\end{equation}}
\newcommand{\bea}{\begin{eqnarray}}
\newcommand{\eea}{\end{eqnarray}}
\newcommand{\ul}{\underline}
\newcommand{\ol}{\overline}
\newcommand{\ti}{\tilde}

\newcommand{\spr}[2]{\langle #1 , #2 \rangle}

\newcommand{\I}{\mathrm{i}}

\newcommand{\re}{\mathrm{Re}}
\newcommand{\im}{\mathrm{Im}}

\newcommand{\res}{\mathrm{Res}}

\newcommand{\bpsi}{\bar{\psi}}
\newcommand{\ulz}{\ul{z}}
\newcommand{\hmu}{\hat{\mu}}

\newcommand{\dimuz}{\di_{\ul{\hat{\mu}}}}
\newcommand{\vrc}{\ul{\Xi}_{p_0}}
\newcommand{\hvrc}{\ul{\hat{\Xi}}_{p_0}}
\newcommand{\di}{\mathcal{D}}

\newcommand{\Amap}{\ul{A}_{p_0}}
\newcommand{\amap}{\ul{\alpha}_{p_0}}
\newcommand{\hAmap}{\ul{\hat{A}}_{p_0}}
\newcommand{\hamap}{\ul{\hat{\alpha}}_{p_0}}

\newcommand{\Rg}[1]{R_{2g+2}^{1/2}(#1)}

\newcommand{\vprod}[2]{\!\!\!\!\begin{array}{c} \mbox{\raisebox{-0.5ex}[0.5ex]
{$\scriptstyle #2 $}} \\ \displaystyle \hspace*{1.1ex}\prod{}^* \\
\mbox{\raisebox{0.6ex}[-0.6ex]{$ \scriptstyle  #1 $}} \end{array}}
\newcommand{\vsum}[2]{\!\!\!\!\begin{array}{c} \mbox{\raisebox{-0.5ex}[0.5ex]
{$\scriptstyle #2 $}} \\ \displaystyle \hspace*{1.1ex}\sum{}^* \\
\mbox{\raisebox{0.6ex}[-0.6ex]{$ \scriptstyle  #1 $}} \end{array}}

\newcommand{\eps}{\varepsilon}

\newcommand{\sig}{\sigma}
\newcommand{\lam}{\lambda}
\newcommand{\gam}{\gamma}
\newcommand{\om}{\omega}

\numberwithin{equation}{section}

\begin{document}

\title[Scattering Theory for Steplike Quasi-Periodic Background]
{Scattering Theory for
Jacobi Operators with Steplike Quasi-Periodic Background}

\author[I. Egorova]{Iryna Egorova}
\address{Kharkiv National University\\ 47,Lenin ave\\ 61164 Kharkiv\\ Ukraine}
\email{\href{mailto:egorova@ilt.kharkov.ua}{egorova@ilt.kharkov.ua}}

\author[J. Michor]{Johanna Michor}
\address{Faculty of Mathematics\\
Nordbergstrasse 15\\ 1090 Wien\\ Austria\\ and International Erwin Schr\"odinger
Institute for Mathematical Physics, Boltzmanngasse 9\\ 1090 Wien\\ Austria}
\email{\href{mailto:Johanna.Michor@esi.ac.at}{Johanna.Michor@esi.ac.at}}
\urladdr{\href{http://www.mat.univie.ac.at/~jmichor/}{http://www.mat.univie.ac.at/\~{}jmichor/}}

\author[G. Teschl]{Gerald Teschl}
\address{Faculty of Mathematics\\
Nordbergstrasse 15\\ 1090 Wien\\ Austria\\ and International Erwin Schr\"odinger
Institute for Mathematical Physics, Boltzmanngasse 9\\ 1090 Wien\\ Austria}
\email{\href{mailto:Gerald.Teschl@univie.ac.at}{Gerald.Teschl@univie.ac.at}}
\urladdr{\href{http://www.mat.univie.ac.at/~gerald/}{http://www.mat.univie.ac.at/\~{}gerald/}}

\thanks{To appear in Inverse Problems.}
\thanks{Work supported by the Austrian Science Fund (FWF) under Grant
No.\ P17762 and INTAS Research Network NeCCA 03-51-6637.}


\begin{abstract}
We develop direct and inverse scattering theory for Jacobi operators
with steplike quasi-periodic finite-gap background in the same isospectral
class. We derive the corresponding Gel'fand-Levitan-Marchenko equation
and find minimal scattering data which determine the perturbed operator
uniquely. In addition, we show how the transmission coefficients can be
reconstructed from the eigenvalues and one of the reflection coefficients.
\end{abstract}

\maketitle

\section{Introduction}

In this paper we consider direct and inverse scattering theory for Jacobi operators with
steplike quasi-periodic finite-gap background in the same isospectral class using the
Marchenko \cite{mar} approach.

Scattering theory for Jacobi operators with a constant background is a classical
topic first developed on an informal level by Case in \cite{dinv4}.
The first rigorous results were established by Guseinov
\cite{gu} with further extensions by Teschl \cite{tivp}, \cite{tjac}.
The case of a quasi-periodic finite-gap background was recently investigated by
us in \cite{emtqps} (see also \cite{voyu} respectively \cite{emtist}, \cite{mtqptr} for applications
to the Toda hierarchy).

Our motivation for the investigation of a steplike situation is twofold. First of all, 
steplike potentials are a simple model in quantum mechanics which have
attracted renewed interest due to their possible applications in mesoscopic solid state
structures. We refer for example to \cite{cvp} where the discrete one-dimensional
Schr\"odinger equation is used as a simple one-band tight binding model to explain
some of the essential qualitative properties of the Wannier-Stark ladders. The
interested reader should also consult \cite{gnp} or \cite{poe} and the references therein
for further information.

Our second motivation is the study of solitons on (quasi-)periodic backgrounds.
While solitons on quasi-periodic backgrounds are well investigated objects,
not much about their stability is known. In fact, as pointed out only recently in \cite{kateptr},
the general believe that the stability problem for solitons on quasi-periodic backgrounds
is similar to the one for solitons on a constant background is wrong (for a detailed
analysis using Riemann-Hilbert techniques see \cite{katept}). This is related
to the fact that solitons on quasi-periodic backgrounds give rise to different spatial asymptotics
which naturally leads  to the type of operators studied here.
Hence our results form the basis for an investigation of solitons on quasi-periodic backgrounds
via the inverse scattering transform. For further details we refer to \cite{emtsr}.

Though the case with different but constant background is well understood
by now (see \cite{eg} and \cite{bdmek}, \cite{dkkz}, \cite{vdo} for applications to the Toda
lattice), the case with different (quasi\mbox{-)}periodic background still is in its infancy.
First results for the case of two period two operators with a special choice for
the respective spectra have been obtained in \cite{baeg}. However, to the best of our
knowledge the general problem is still open (even for the
case of one-dimensional Schr\"odinger operators -- see \cite{gnp} and the
references therein for a recent account). Our aim is to fill this gap for the case
of isospectral background operators.

After recalling some necessary facts on algebro-geometric quasi-periodic finite-gap
operators in Section~\ref{secQP}, we construct the transformation operators and 
investigate the properties of the scattering data in Section~\ref{secSD}. 
In particular, we show
how both transmission coefficients can be reconstructed from the eigenvalues and
the corresponding reflection coefficients. In Section~\ref{secGLM} we derive the 
Gel'fand-Levitan-Marchenko equation and show that it uniquely determines the
operator. In addition, we formulate
necessary conditions for the scattering data to uniquely determine our Jacobi operator.
Our final Section~\ref{secINV} shows that our necessary conditions for the scattering
data are also sufficient.

\section{Quasi-periodic finite-gap operators}
\label{secQP}

As a preparation for our next section we first need to recall some facts on
quasi-periodic finite-gap Jacobi operators which contain all
periodic operators as a special case. We refer to \cite[Chapter~9]{tjac}.

Let $\M$ be the Riemann surface associated with the following function
\begin{equation}
\Rg{z}= -\prod_{j=0}^{2g+1} \sqrt{z-E_j}, \qquad
E_0 < E_1 < \cdots < E_{2g+1},
\end{equation}
where $g\in \N$ and $\sqrt{.}$ is the standard root with branch cut along $(-\infty,0)$. 
$\M$ is a compact, hyperelliptic Riemann surface of genus $g$.
A point on $\M$ is denoted by 
$p = (z, \pm \Rg{z}) = (z, \pm)$, $z \in \C$, or $p = \infty_\pm$, and
the projection onto $\C \cup \{\infty\}$ by $\pi(p) = z$. The sets 
$
\Pi_\pm = \{ (z, \pm \Rg{z}) \mid z \in \C\backslash
\bigcup_{j=0}^g[E_{2j}, E_{2j+1}]\} \subset \M
$
are called upper, lower sheet, respectively.

Let $\{a_j, b_j\}_{j=1}^g$ be loops on the Riemann surface $\M$ representing the
canonical generators of the fundamental group $\pi_1(\M)$. We require
$a_j$ to surround the points $E_{2j-1}$, $E_{2j}$ (thereby changing sheets
twice) and $b_j$ to surround $E_0$, $E_{2j-1}$ counter-clockwise on the
upper sheet, with pairwise intersection indices given by
\begin{equation}
a_j \circ a_k= b_j \circ b_k = 0, \qquad a_j \circ b_k = \delta_{jk},
\qquad 1 \leq j, k \leq g.
\end{equation}
The corresponding canonical basis $\{\zeta_j\}_{j=1}^g$ for the space of
holomorphic differentials can be constructed by
\begin{equation}
\underline{\zeta} = \sum_{j=1}^g \underline{c}(j)  
\frac{\pi^{j-1}d\pi}{R_{2g+2}^{1/2}},
\end{equation}
where the constants $\underline{c}(.)$ are given by
\[
c_j(k) = C_{jk}^{-1}, \qquad 
C_{jk} = \int_{a_k} \frac{\pi^{j-1}d\pi}{R_{2g+2}^{1/2}} =
2 \int_{E_{2k-1}}^{E_{2k}} \frac{z^{j-1}dz}{\Rg{z}} \in
\R.
\]
The differentials fulfill
\begin{equation}
\int_{a_j} \zeta_k = \delta_{j,k}, \qquad \int_{b_j} \zeta_k = \tau_{j,k}, 
\qquad \tau_{j,k} = \tau_{k, j}, \qquad 1 \leq j, k \leq g.
\end{equation}
Finally we will need $\om_{p q}$, the normalized Abelian differential of the third kind with poles
at $p$ and $q$, that is, $\om_{p q}$ has vanishing $a$-periods and first order poles with residues
$+1$, $-1$ at $p$, $q$, respectively.

Now pick $g$ numbers (the Dirichlet eigenvalues)
\be
(\hat{\mu}_j)_{j=1}^g = (\mu_j, \sigma_j)_{j=1}^g
\ee
whose projections lie in the spectral gaps, that is, $\mu_j\in[E_{2j-1},E_{2j}]$.
Associated with these numbers is the divisor $\dimuz$ which
is one at the points $\hat{\mu}_j$  and zero else. Using this divisor we
introduce
\begin{align} \nn
\ulz(p,n) &= \hAmap(p) - \hamap(\dimuz) - n\ul{\hat A}_{\infty_-}(\infty_+)
- \hvrc \in \C^g, \\ \label{ulz}
\ulz(n) &= \ulz(\infty_+,n),
\end{align}
where $\vrc$ is the vector of Riemann constants and $\Amap$ ($\amap$)
is Abel's map (for divisors). The hat indicates that we
regard it as a (single-valued) map from $\hat{\M}$ (the fundamental polygon
associated with $\M$) to $\C^g$.
We recall that the function $\theta(\ulz(p,n))$ has precisely $g$ zeros
$\hmu_j(n)$ (with $\hmu_j(0)=\hmu_j$), where $\theta(\ul{z})$ is the
Riemann theta function of $\M$.

With this notation our quasi-periodic finite-gap operator is given by
\be
H_q f(n) = a_q(n) f(n+1) + a_q(n-1) f(n-1) + b_q(n) f(n),
\ee
where
\begin{align} \nn
a_q(n)^2 &= \ti{a}^2 \frac{\theta(\ulz(n+1)) \theta(\ulz(n-1))}{\theta(
\ulz(n))^2},\\ \label{imfab}
b_q(n) &= \tilde{b} + \sum_{j=1}^g c_j(g)
\frac{\partial}{\partial w_j} \ln\Big(\frac{\theta(\ul{w} +
\ulz(n)) }{\theta(\ul{w} + \ulz(n-1))}\Big) \Big|_{\ul{w}=0}.
\end{align}
The constants $\ti{a}$, $\tilde{b}$, $c_j(g)$ depend only on the Riemann surface
(see \cite[Section~9.2]{tjac}).

Introduce the Baker-Akhiezer function
\begin{align} \nn \label{BAfthetarep}
\psi_q(p,n) &= C(n,0) \frac{\theta (\ulz(p,n))}{\theta(\ulz (p,0))}
\exp \Big( n \int_{p_0}^p \om_{\infty_+,\infty_-}
\Big),
\end{align}
where $C(n,0)$ is real-valued,
\begin{equation}
C(n,0)^2 = \frac{ \theta(\ulz(0)) \theta(\ulz(-1))}
{\theta (\ulz (n))\theta (\ulz (n-1))},
\end{equation}
and $\om_{\infty_+,\infty_-}$ is the Abelian differential of the third kind with poles
at $\infty_+$ respectively $\infty_-$. The two branches $\psi_{q,\pm}(z,n)= \psi_q(p,n)$,
$p=(z,\pm)$, of the Baker-Akhiezer function are solutions (in the weak sense) of
\begin{align}  \nonumber
H_q \psi_{q,\pm}(z,n) &= z \psi_{q,\pm}(z,n)
\end{align}
which are linearly independent away from the branch points. In fact, their Wronskian
is given by
\begin{align}  \nonumber
W_q(\psi_{q,-}(z), \psi_{q,+}(z)) &= a_q(n) \left(\psi_{q,-}(z,n), \psi_{q,+}(z,n+1) -
\psi_{q,-}(z,n+1), \psi_{q,+}(z,n)\right)\\
&= \frac{R^{1/2}_{2g+2}(z)}{\prod_{j=1}^g (z-\mu_j)}.
\end{align}
It is well-known that the spectrum of $H_q$ consists of $g+1$ bands
\begin{equation}
\sig(H_q) = \bigcup_{j=0}^g [E_{2j},E_{2j+1}].
\end{equation}
Note that $\psi_{q,\pm}(z,n)$ are discontinuous along the spectrum and we set
$\psi_{q,\pm}(\lam,n)= \lim_{\eps\downarrow0} \psi_{q,\pm}(\lam+\I\eps,n)$,
$\lam\in\sig(H_q)$. This implies $\lim_{\eps\downarrow0} \psi_{q,\pm}(\lam-\I\eps,n)=
\ol{\psi_{q,\pm}(\lam,n)} = \psi_{q,\mp}(\lam,n)$, $\lam\in\sig(H_q)$.
For further information and proofs we refer to \cite[Chapter~9]{tjac}.

\section{Scattering data}
\label{secSD}

Consider  two quasi-periodic
finite-gap operators $H_q^\pm$ associated with the sequences $a_q^\pm$, $b_q^\pm$
in the same isospectral class,
\be
\sig(H_q^+) = \sig(H_q^-) \equiv \Sigma = \bigcup_{j=0}^g [E_{2j},E_{2j+1}],
\ee
but with possibly different Dirichlet data $\{\hat{\mu}_j^\pm\}_{j=1}^g$. We will
add $\pm$ as a superscript to all data introduced in Section~\ref{secQP} to
distinguish between the corresponding data of $H_q^+$ and $H_q^-$.
To avoid excessive  sub/superscripts we abbreviate
\be
\psi_q^\pm(z, n)= \psi_{q,\pm}^\pm(z, n) \quad\mbox{and}\quad
\bpsi_q^\pm(z, n)= \psi_{q,\mp}^\pm(z, n),
\ee
that is, $\psi_q^\pm(z, n)$ is the solution of $H_q^\pm$ decaying near $\pm\infty$ and
$\bpsi_q^\pm(z, n)$ is the solution of $H_q^\pm$ decaying near $\mp\infty$.
Note that for $\lam\in\Sigma$ we have $\bpsi_q^\pm(\lam,n)= \ol{\psi_q^\pm(\lam,n)}$.

In addition, we split the set of Dirichlet eigenvalues into three parts:
\begin{align} \nn
M^\pm &= \{ \mu_j^\pm \,|\,
\mu_j^\pm \in \R\backslash\Sigma \mbox{ is a pole of } \psi_q^\pm(z,1)\},\\ \nn
\bar{M}^\pm &= \{ \mu_j^\pm \,|\,
\mu_j^\pm  \in \R\backslash\Sigma \mbox{ is a pole of } \bpsi_q^\pm(z,1)\},\\
\tilde{M}^\pm &= \{ \mu_j^\pm \,|\, \mu_j^\pm \in \partial \Sigma \},
\end{align}
and introduce
\be
\hat{\psi}_q^\pm(z,n) = \bigg(\prod_{\mu\in M^\pm} (z-\mu) \bigg) \psi_q^\pm(z,n),
\ee
which are holomorphic for all $z\in\C\backslash\Sigma$.

Let $a(n)$, $b(n)$ be sequences satisfying
\be \label{hypo}
\sum_{n = 0}^{\pm \infty} |n| \Big(|a(n) - a_q^\pm(n)| + |b(n) - b_q^\pm(n)| \Big)
< \infty
\ee
and denote the corresponding operator by $H$.

The special case $H_q^- = H_q^+$ has been exhaustively studied in \cite{emtqps}
(see also \cite{voyu})
and several results are straightforward generalizations. In such situations we will
simply refer to \cite{emtqps} and only point out possible differences.
In fact, we begin by recalling two theorems from \cite{emtqps}:

\begin{theorem} \label{thmjost}
Assume (\ref{hypo}). Then there exist weak solutions
$\psi_\pm(z, .)$ of $H \psi = z \psi$ satisfying
\begin{equation}                         \label{perturbed sol}
\lim_{n \rightarrow \pm \infty}
|w(z)^{\mp n} (\psi_\pm(z, n) - \psi_q^\pm(z, n))| = 0,
\end{equation}
where $w(z)= \exp(\int_{p_0}^{(z,+)} \om_{\infty_+,\infty_-})$ is the quasi-momentum map.
Moreover, $\psi_\pm(z, .)$ are continuous 
(resp.\ holomorphic) with respect to $z$ whenever $\psi_q^\pm(z, .)$ are and
\be                 \label{B4jost}
\psi_\pm(z,n) =   \frac{z^{\mp n}}{A_\pm(n)} \Big(\vprod{j=0}{n-1}a_q^\pm(j)\Big)^{\pm 1} 
\Big(1 + \Big(B_\pm(n) \pm \vsum{j=1}{n} b_q^\pm(j- {\scriptstyle{0 \atop 1}}) \Big)\frac{1}{z}
+ O(\frac{1}{z^2}) \Big),
\ee
where
\bea \nonumber
A_+(n) &=& \prod_{j=n}^\infty \frac{a(j)}{a_q^+(j)}, \qquad
B_+(n)= \sum_{m=n+1}^\infty (b_q^+(m)-b(m)), \\
A_-(n) &=& \prod_{j=- \infty}^{n-1} \frac{a(j)}{a_q^-(j)}, \qquad
B_-(n) = \sum_{m=-\infty}^{n-1} (b_q^-(m)-b(m)).
\eea
\end{theorem}

\noindent
Defining $\hat\psi_\pm(z,n)$ using $\hat\psi_q^\pm(z,n)$ (instead of $\psi_q^\pm(z,n)$), we
have that $\hat\psi_\pm(z,n)$ are holomorphic in $\C\backslash\Sigma$.

\begin{theorem} 
Assume (\ref{hypo}). Then we have $\sig_{ess}(H)=\Sigma$, the
point spectrum of $H$ is finite and confined to the spectral gaps of
$H_q^\pm$, that is, $\sig_p(H) = \{ \rho_j\}_{j=1}^q \subset \R\backslash\Sigma$.
Furthermore, the essential spectrum of $H$ is purely absolutely continuous.
\end{theorem}

\noindent
\noindent
It will be convenient to identify $\C\backslash\Sigma$ with $\Pi_+$ and
regard the functions $\psi_\pm$, $\psi_q^\pm$ as functions on $\Pi_+$.
Furthermore, we extend them by continuity to $\partial\Pi_+$, that is,
$\psi_+((\lam,\pm),n) =  \lim_{\eps\downarrow0} \psi_+(\lam\pm\I\eps,n)$,
etc, for $\lam\in\Sigma$. Consequently we have $\psi_q^\pm(p^*,n) =
\ol{\psi_q^\pm(p,n)}$ and $\psi_\pm(p^*,n) =
\ol{\psi_\pm(p,n)}$ for $p\in\partial\Pi_+$.

Note: Unfortunately this disagrees with our previous definition of $\psi_q^-$ which
was originally defined as the branch on $\Pi_-$. Moreover, one should emphasize
that $\psi_\pm(p,n)$ does not have a continuation to $\Pi_-$ in general.

Using the fact that $\psi_q^\pm(p,n)$ form an orthonormal basis for
$L^2(\partial \Pi_+,d\omega^\pm)$, where
\be \label{domega}
d\omega^\pm = \frac{\prod_{j=1}^g(\pi-\mu_j^\pm)}{R^{1/2}_{2g+2}} d\pi,
\ee
we can define  
\begin{align} \nn
K_\pm(n,m) &= \int_{\partial\Pi_+} \psi_\pm(p,n) \ol{\psi_q^\pm(p, m)}
d\omega^\pm \\ 
&= 2 \re \int_\Sigma \psi_\pm(\lam,n) \ol{\psi_q^\pm(\lam, m)}
d\omega^\pm.
\end{align}

Then $K_\pm(n,m)$ satisfy (\cite{emtqps}):

\begin{lemma}       \label{le:jost kernel q}
Assume (\ref{hypo}).
The Jost solutions $\psi_\pm(z,n)$ can be represented as  
\be
\psi_\pm(z,n) = \sum_{m=n}^{\pm \infty} K_\pm(n,m)
\psi_q^\pm(z, m),
\ee
where the kernels $K_\pm(n, .)$ satisfy
$K_\pm(n,m) = 0$ for $\pm m < \pm n$ and  
\be                      \label{estimate K q}
|K_\pm(n,m)| \leq C_\pm(n) \sum_{j=[\frac{m+n}{2}] \pm 1}^{\pm \infty}
\Big(|a(j)-a_q^\pm(j)| + |b(j)-b_q^\pm(j)|\Big), \quad \pm m > \pm n>0.
\ee
The functions $C_\pm(n) > 0$ decrease as $n \rightarrow \pm \infty$ and depend only on 
$H_q^\pm$ and the value of the sums in (\ref{hypo}).
\end{lemma}

Associated with $K_\pm(n,m)$ is the operator 
\be
(\mathcal{K_\pm} f)(n) = \sum_{m=n}^{\pm \infty} K_\pm(n,m) f(m), 
\qquad  f \in \ell_\pm^\infty(\Z, \C),
\ee
which acts as a transformation operator for
the pair $\tau$, $\tau_q^\pm$. 

\begin{theorem} \label{th:transf}
Let $\tau_q^\pm$ and $\tau$ be the quasi-periodic and perturbed 
Jacobi difference expressions, respectively. Then \label{th:transf op q}
\be
\tau \mathcal{K_\pm}f = \mathcal{K_\pm} \tau_q^\pm f, \qquad f \in
\ell_\pm^\infty(\Z, \C).
\ee
Furthermore, for $n \in \Z$ we have
\begin{align}              \label{a a_q}
\frac{a(n)}{a_q^+(n)} &= \frac{K_+(n+1, n+1)}{K_+(n,n)}\\   \nonumber
\frac{a(n)}{a_q^-(n)} &= \frac{K_-(n,n)}{K_-(n+1, n+1)},  \\          \nonumber
b(n) - b_q^+(n) &= a_q^+(n)\frac{K_+(n, n+1)}{K_+(n,n)}
-a_q^+(n-1)\frac{K_+(n-1, n)}{K_+(n-1,n-1)}  \\           \nonumber
b(n) - b_q^-(n) &= a_q^-(n-1)\frac{K_-(n, n-1)}{K_-(n,n)}
-a_q^-(n)\frac{K_-(n+1, n)}{K_-(n+1,n+1)}. 
\end{align}
\end{theorem}

\noindent
Next we define the coefficients of the scattering matrix via the scattering relations
\be
\psi_\mp (\lam,n) = \alpha_\pm(\lam) \ol{\psi_\pm(\lam,n)} 
+ \beta_\pm(\lam) \psi_\pm(\lam,n), \qquad \lam \in \Sigma,
\ee
where
\begin{align} \label{alpha quasi}
\alpha_\pm(\lam) &= \frac{W(\psi_\pm(\lam), \psi_\mp(\lam))}
{W(\psi_\pm(\lam), \ol{\psi_\pm(\lam)})} 
= \frac{\prod_{j=1}^g (\lam-\mu_j^\pm)}{\Rg{\lam}} W(\psi_-(\lam), \psi_+(\lam)), \\ \nn
\beta_\pm(\lam) &= \frac{W(\psi_\mp(\lam), \ol{\psi_\pm(\lam)})}
{W(\psi_\pm(\lam), \ol{\psi_\pm(\lam)})} 
= \mp \frac{\prod_{j=1}^g (z-\mu_j^\pm)}{\Rg{z}}W(\psi_\mp(\lam), \ol{\psi_\pm(\lam)}), 
\end{align}
and $W_n(f,g)=a(n)(f(n)g(n+1)-f(n+1)g(n))$ denotes the Wronskian.
Transmission $T_\pm(\lam)$ and reflection $R_\pm(\lam)$ 
coefficients are then defined by
\be
T_\pm(\lam) = \alpha_\pm^{-1}(\lam), \qquad
R_\pm(\lam) = \frac{\beta_\pm(\lam)}{\alpha_\pm(\lam)} = 
\frac{W(\psi_\mp(\lam), \ol{\psi_\pm(\lam)})}{W(\psi_\pm(\lam), \psi_\mp(\lam))}.
\ee

Using
\be \label{T}
T_\pm(z) =  \frac{\Rg{z}}{\prod_{j=1}^g (z-\mu_j^\pm)} \frac{1}{W(\psi_-(z), \psi_+(z))}
\ee
we see that $T_\pm(z)$ admit a meromorphic extension to $\C\backslash\Sigma$.
Since $W(\hat\psi_-(z),\hat\psi_+(z))$ is holomorphic in $\C\backslash\Sigma$
with simple zeros at the eigenvalues $\rho_k$ (see \cite[Section 2.2]{tjac})
we obtain the following behavior:
\be
T_\pm(z) =
\frac{\prod_{\mu_j^\mp\in M^\mp}(z-\mu_j^\mp)}{\prod_{\mu_j^\pm\in\bar{M}^\pm}(z-\mu_j^\pm)}
\frac{D_\pm(z)}{\prod_{k=1}^q (z-\rho_k)},
\ee
where $D_\pm(z)$ are holomorphic and do not vanish in $\C\backslash\Sigma$.
That is, in general the transmission coefficients have simple poles at the eigenvalues
$\rho_j$ of $H$. In addition, there are simple poles at $\mu_j^\pm\in \bar{M}^\pm$ and simple
zeros at $\mu_j^\mp\in M^\mp$. A pole at $\mu_j^\pm$
could cancel with a zero at $\mu_j^\mp$ or could give a second order pole if
$\mu_j^\pm=\rho_k$.

The Pl\"ucker identity (c.f.\ \cite[(2.169)]{tjac}) implies
\be
\alpha_+(\lam) \ol{\alpha_-(\lam)} = 1 - \beta_+(\lam) \beta_-(\lam),
\ee
and using 
\be
\alpha_+(\lam) = \alpha_-(\lam) \prod_{j=1}^g\frac{\lam - \mu_j^+}{\lam - \mu_j^-}, \qquad
\ol{\beta_+}(\lam) = - \beta_-(\lam) \prod_{j=1}^g\frac{\lam - \mu_j^+}{\lam - \mu_j^-},
\ee
we obtain
\be
|\alpha_\pm(\lam)|^2 = \prod_{j=1}^g\frac{\lam - \mu_j^\pm}{\lam - \mu_j^\mp} + |\beta_\pm(\lam)|^2.
\ee
The norming constants $\gam_{\pm, j}$ corresponding to $\rho_j \in \sigma_p(H)$ 
are given by
\be \label{norming}
\gam_{\pm, j}^{-1}=\sum_{n \in \Z}|\hat\psi_\pm(\rho_j, n)|^2, \qquad 1 \leq j \leq q.
\ee
Moreover, we set $\hat\psi_\pm(\rho_j, .)= c_j^\pm \hat\psi_\mp(\rho_j, .)$ with $c_j^+ c_j^-=1$.

\begin{lemma}
The coefficients $T_\pm(\lam)$, $R_\pm(\lam)$ are continuous for
$\lam \in \Sigma$ except at possibly the band edges $E_j$, and fulfill
\begin{align}  \label{TR1}
T_+(\lam) \ol{T_-(\lam)} + |R_\pm(\lam)|^2 &= 1, \qquad \lam \in \Sigma,\\
T_\pm(\lam) \ol{R_\pm(\lam)} + \ol{T_\pm(\lam)} R_\mp(\lam) &= 0, \qquad \lam \in \Sigma.
\end{align}
In particular,
\be \label{TR3}
|T_\pm(\lam)|^2  \prod_{j=1}^g\frac{\lam - \mu_j^\pm}{\lam - \mu_j^\mp} 
+ |R_\pm(\lam)|^2 =1,
\ee
and hence $|R_\pm(\lam)|^2\le 1$ with equality only possibly at the band edges $E_j$,
where
\be \label{T E}
\begin{array}{r@{\qquad}l}
\lim\limits_{z \rightarrow E} \Rg{z}\frac{R_\pm(z) + 1}{T_\pm(z)} =0,
& E\neq \{\mu_j^+, \mu_j^-\},\\
\lim\limits_{z \rightarrow E} R_{2g+2}(z) \frac{R_\pm(z) + 1}{T_\pm(z)} =0,
& E = \mu_j^\mp \neq \mu_j^\pm, \\
\lim\limits_{z \rightarrow E} \frac{R_\pm(z) - 1}{T_\pm(z)} =0,
& E = \mu_j^\pm \neq \mu_j^\mp, \\
\lim\limits_{z \rightarrow E} \Rg{z} \frac{R_\pm(z) - 1}{T_\pm(z)} =0,
& E = \mu_j^+ = \mu_j^-.
\end{array}
\ee
The transmission coefficients $T_\pm(\lam)$ have a meromorphic continuation 
with
\be \label{res T}
\res_{z=\rho_k}
\frac{\prod_{\mu_j^\mp\in M^\mp}(z - \mu_j^\mp)}{\prod_{\mu_j^\pm\in \bar{M}^\pm \cup \tilde{M}^\pm}(z - \mu_j^\pm)}
T_\pm(z) = -\sqrt{\gam_{+,k}\gam_{-,k}} \Rg{\rho_k}.
\ee
In addition, $T_\pm(z) \in \R$ as $z \in \R \backslash \Sigma$ and
\be \label{T infty}
T_\pm(\infty) = \frac{1}{K_+(0,0) K_-(0,0)}.
\ee
\end{lemma}

\begin{proof}
To prove (\ref{T E}), recall first that if $E \notin \ti M^\pm$, 
then $\psi_\pm(\lam,n)- \ol{\psi_\pm(\lam,n)} \to 0$ as $\lam \to E$, 
and $\psi_\pm(E,n)$ are bounded.
If $E \in \ti M^\pm$, then $R_{2g+2}^{1/2}(\lam)\psi_\pm(\lam,n)$ and
$\psi_\pm(\lam,n)+ \ol{\psi_\pm(\lam,n)}$ are bounded at $E$.
 
Thus using definition (\ref{alpha quasi}), we have as $E \notin \ti M^- \cup \ti M^+$
\[
R^{1/2}_{2g+2}(\lambda)  \frac{R_\pm(\lambda) + 1}{T_\pm(\lambda)} =
\prod_{j=1}^g (\lambda - \mu_j^\pm) \big(W(\psi_-(\lambda), \psi_+(\lambda)) \mp
W(\psi_\mp(\lambda), \ol{\psi_\pm(\lambda)})\big) \to 0.
\] 
If $E \in \ti M^\pm \backslash \ti M^\mp$, then
\[
\frac{R_\pm(\lambda) - 1}{T_\pm(\lambda)} = \pm
\frac{\prod_{j=1}^g (\lambda - \mu_j^\pm)}{R^{1/2}_{2g+2}(\lambda)} 
W(\psi_\pm(\lambda) + \ol{\psi_\pm(\lambda)}, \psi_\mp(\lambda)) \to 0.
\]
In the same manner one proves the remaining two equalities in (\ref{T E}).

By \cite[(2.33)]{tjac},
\be \label{Wprime}
W^\prime(\hat\psi_-(z), \hat\psi_+(z)) \Big|_{z=\rho_j} = \frac{-1}{\sqrt{\gam_{-,j} \gam_{+,j}}},
\ee
which proves (\ref{res T}). Finally, (\ref{T infty}) follows from (\ref{B4jost}).
\end{proof}

\noindent
Observe that (\ref{TR1}) implies $|R_-(\lam)|=|R_+(\lam)|$.

The sets  
\be
S_\pm(H) = \{R_\pm(\lam), \lam \in \Sigma; \, (\rho_j, \gam_{\pm, j}), 
1 \leq j \leq q\}
\ee
are called left/right scattering data for $H$.

\begin{theorem}
The transmission coefficients can be reconstructed from one reflection
coefficient and the eigenvalues via
\begin{align}\nn
T_+(z) =&
\sqrt{\frac{\theta(\ulz^+(\infty_+,0))\theta(\ulz^+(\infty_-,0))}{\theta(\ulz^-(\infty_+,0)) \theta(\ulz^-(\infty_-,0))}}
\frac{\theta(\ulz^-(\hat z^*,0))}{\theta(\ulz^+(\hat z^*,0))} \times \\  \label{rectp}
 & \times \left( \prod_{j=1}^q \exp\left(-\int_{E(\rho_j)}^{\hat\rho_j} \!\!\! \om_{\hat z \hat z^*}\right) \right)
\exp \left(\frac{1}{2\pi \I} \int_{\Sigma}  
\ln(1 - |R_\pm|^2) \om_{\hat z \hat z^*} \right),\\ \nn
T_-(z) =&
\sqrt{\frac{\theta(\ulz^-(\infty_+,0)) \theta(\ulz^-(\infty_-,0))}{\theta(\ulz^+(\infty_+,0))\theta(\ulz^+(\infty_-,0))}}
\frac{\theta(\ulz^+(\hat z,0))}{\theta(\ulz^-(\hat z,0))} \times \\  \label{rectm}
 & \times \left( \prod_{j=1}^q \exp\left(-\int_{E(\rho_j)}^{\hat\rho_j} \!\!\! \om_{\hat z \hat z^*}\right) \right)
\exp \left(\frac{1}{2\pi \I} \int_{\Sigma}  
\ln(1 - |R_\pm|^2) \om_{\hat z \hat z^*} \right),
\end{align}
where $\hat z =(z,+)$, the integral over $\Sigma$ is taken on the upper sheet, and
$E(\rho)$ is $E_0$ if $\rho<E_0$, either $E_{2j-1}$ or $E_{2j}$ if
$\rho\in(E_{2j-1},E_{2j})$, $1\le j \le g$, and $E_{2g+1}$ if $\rho>E_{2g+1}$.

Furthermore, the phase shift between $H_q^+$ and $H_q^-$ can also
be computed from one reflection coefficient and the eigenvalues
\be \label{phaseshift}
\amap(\di_{\ul{\hat{\mu}}^+}) - \amap(\di_{\ul{\hat{\mu}}^-})
= \sum_{j=1}^q \int_{\hat\rho_j^*}^{\hat\rho_j} \ul\zeta -
\frac{1}{2\pi\I} \int_{\partial \Pi_+}  \!\!\! \ln(1 - |R_\pm|^2) \ul\zeta.
\ee
\end{theorem}

\begin{proof}
First of all note that $\im \int \om_{pp^*}$ is the Green's function of $\Pi_+$ and that
\be
B(z,\rho)= \exp\left(\int_{E(\rho)}^{\hat\rho} \!\!\! \om_{\hat z \hat z^*}\right)
\ee
is the Blaschke factor (see \cite{tistalg}). That is, $B(z,\rho)$ is a multivalued holomorphic function,
which vanishes at $z=\rho$ and satisfies $|B(\lam,\rho)|=1$ for $\lam\in\Sigma$.

We start by considering the multivalued function
\be
t_+(z) =
\sqrt{\frac{\theta(\ulz^-(\infty_+,0)) \theta(\ulz^-(\infty_-,0))}{\theta(\ulz^+(\infty_+,0))\theta(\ulz^+(\infty_-,0))}}
\frac{\theta(\ulz^+(\hat z^*,0))}{\theta(\ulz^-(\hat z^*,0))}
\left( \prod_{j=1}^q B(z,\rho_j) \right) T_+(z)
\ee
which has neither zeros nor poles on $\Pi_+$. 
Using the same argument as in \cite[eq. (9.27)]{tjac} one verifies that
$\theta(\ulz^+(\infty_\pm,0)) / \theta(\ulz^-(\infty_\pm,0))$ is positive.
The absolute value of $t_+(z)$ is single-valued
and, using $\ol{\ulz^\pm(\hat\lam,n)} = \ulz^\pm(\hat\lam^*,n) \mod\Z^g$ for 
$\hat\lam =(\lam,+)$ with $\lam\in\Sigma$,
one obtains
\begin{align} \nn
|t_+(\lam)|^2 &=
\frac{\theta(\ulz^-(\infty_+,0)) \theta(\ulz^-(\infty_-,0))}{\theta(\ulz^+(\infty_+,0))\theta(\ulz^+(\infty_-,0))}
\frac{\theta(\ulz^+(\hat\lam^*,0))}{\theta(\ulz^-(\hat\lam^*,0))}
\frac{\theta(\ulz^+(\hat\lam,0))}{\theta(\ulz^-(\hat\lam,0))} 
|T_+(\lam)|^2\\
& = \frac{\prod (\lam-\mu_j^+)}{\prod (\lam-\mu_j^-)} |T_+(\lam)|^2.
\end{align}
Here the last equality follows since the ratio of theta functions extends to a
(single-valued) meromorphic function on $\C$ which is one at $\infty$.
Hence $|t_+(\lam)|^2 = 1 - |R_\pm(\lam)|^2$ for $\lam\in\Sigma$ by (\ref{TR3}) and
$\ln|t_+(z)|$ can be reconstructed from its boundary values using Green's function:
\be
|t_+(z)| = \exp \left(\re\frac{1}{2\pi \I} \int_{\Sigma}  
\ln(1 - |R_\pm|^2) \om_{\hat z \hat z^*} \right).
\ee
So (\ref{rectp}) holds at least when taking absolute values. But since both sides are
meromorphic, they can only differ by a factor of absolute value one. Since
both sides are positive at $\infty_+$, they are equal. The claim for $T_-(z)$
follows analogously.

To show (\ref{phaseshift}) we use that, since $T_\pm$ are single-valued, there is no
jump when we go around $b$-cycles. Hence the jump when going around $b_\ell$ of all
factors must add up to zero, which is just (\ref{phaseshift}).
\end{proof}

\noindent
Combining this result with (\ref{res T}) we obtain:

\begin{corollary} \label{corsd}
One of the scattering data $S_-(H)$ or $S_+(H)$ determines the other.
\end{corollary}

\section{The Gel'fand-Levitan-Marchenko equation}
\label{secGLM}

Finally we want to derive the Gel'fand-Levitan-Marchenko equation and
show that one of the sets $S_-(H)$ or $S_+(H)$  uniquely determines $H$.

\begin{theorem}
The kernel $K_\pm(n,m)$ of the transformation operator satisfies
the Gel'fand-Levitan-Marchenko equation
\be   \label{glm1 q}
K_\pm(n,m) + \sum_{l=n}^{\pm \infty}K_\pm(n,l)F^\pm(l,m) = 
\frac{\delta(n,m)}{K_\pm(n,n)},  \qquad \pm m \geq \pm n,
\ee
where
\begin{align} \nn
F^\pm(m,n) =& 2\re \int_\Sigma R_\pm(\lam) \psi_q^\pm(\lam,m) 
\psi_q^\pm(\lam,n) d\omega^\pm\\ \label{glm2 q}
& + \sum_{j=1}^q \gam_{\pm,j}  \hat\psi_q^\pm(\rho_j,n) \hat\psi_q^\pm(\rho_j,m).
\end{align}
\end{theorem}

\begin{proof}
As already done for $\psi_\pm$, we regard $T_\pm(z)$ as functions on $\Pi_+$
and $R_\pm(\lam)$ as functions on $\partial\Pi_+$ by setting
$R_\pm((\lam,+))= R_\pm(\lam)$ respectively $R_\pm((\lam,-))= \ol{R_\pm(\lam)}$.

Computing the Fourier coefficients of
\be
T_\pm(\lam) \psi_\mp(\lam, n) = R_\pm(\lam) \psi_\pm(\lam, n) +
\ol{\psi_\pm(\lam, n)} 
\ee
as in \cite[Sec.~7]{emtqps} one obtains by (\ref{B4jost}), (\ref{domega}), (\ref{T}), (\ref{T infty}),
(\ref{Wprime}),
$\hat\psi_-(\rho_j, n)= c_j^- \hat\psi_+(\rho_j, n)$, and the residue theorem
\begin{align*}     \nonumber
&\int_{\partial\Pi_+} T_+(p) \psi_-(p, n)\psi_q^+(p, m)d\omega^+ 
= \int_{\partial\Pi_+}
\frac{\psi_-(p, n)\psi_q^+(p, m)}{W(\psi_-(p), \psi_+(p))}d\pi \\ \nonumber
& \quad = \frac{\delta(n,m)}{K_+(n,n)} 
+ \sum_{j=1}^q  \res_{\rho_j} 
\bigg(
\frac{\hat\psi_-(p, n) \hat\psi_q^+(p, m)}{W(\hat\psi_-(p), \hat\psi_+(p))}
\bigg)\\ \nonumber
& \quad = \frac{\delta(n,m)}{K_+(n,n)} 
- \sum_{j=1}^q \gam_{+,j} \hat\psi_+(\rho_j,n) \hat\psi_q^+(\rho_j,m).
\end{align*}
The right hand side follows analogous to \cite{emtqps}. 
\end{proof}

Note that while the scattering data depend on the particular normalization chosen
for $\psi_q^\pm$, the kernel of Gel'fand-Levitan-Marchenko equation is of course independent
of this normalization.

\begin{theorem} \label{thmglm}
For $n \in \Z$, the Gel'fand-Levitan-Marchenko operator
\be
\mathcal{F}^\pm_n: \ell^2 \to \ell^2, \qquad 
\mathcal{F}^\pm_n f(j) = \sum_{l=0}^\infty F^\pm(n \pm l, n \pm j)f(l),
\ee
is Hilbert-Schmidt.
Moreover, $1+\mathcal{F}^\pm_n$ is positive and hence invertible.

In particular, the Gel'fand-Levitan-Marchenko equation (\ref{glm1 q}) 
\be       \label{marchenko 2 q}              
(1 + \mathcal{F}^\pm_n) K_\pm(n, n \pm .) = (K_\pm(n, n))^{-1} \delta_0
\ee
has a unique solution and $S_+(H)$ or $S_-(H)$ uniquely determine $H$.
\end{theorem}

\begin{proof}
The proof can be done as in \cite[Theorem~7.5]{emtqps}, since
by equation (\ref{TR3}) we have $|R_\pm(\lam)|<1$ for 
$\lam \in \Sigma \backslash \partial \Sigma$. 
\end{proof}

To finish the direct scattering step we summarize the properties of
$S_{\pm}(H)$.

\begin{hypo}             \label{hypo scat q}
The scattering data
\[
S_{\pm}(H) = \{R_{\pm}(\lam), \lam \in \Sigma; (\rho_j, \gamma_{\pm, j}), 1 \leq j
\leq q\}
\]
satisfy the following conditions: 

(i) The reflection coefficients $R_{\pm}(\lam)$ are continuous except
possibly at the band edges $E$ and fulfill
$|R_\pm(\lam)|<1$ for $\lam \in \Sigma \backslash \partial \Sigma$.

The Fourier coefficients
\be
\tilde F^\pm(l,m) = 2 \re\int_{\Sigma} R_\pm(\lam) \psi_q^\pm(\lam,l) 
\psi_q^\pm(\lam,m) d\omega^\pm
\ee
satisfy  
\bea         \nonumber                      
&&|\tilde F^{\pm}(n, m)| \leq \sum_{j=n+m}^{\pm \infty} q(j), 
\qquad q(j) \geq 0, \qquad  |j| q(j) \in \ell^1(\Z), \\         \nonumber
&&\sum_{n = n_0}^{\pm \infty}|n| \big|\tilde F^{\pm}(n,n) - 
\tilde F^{\pm}(n \pm 1, n \pm 1)\big| < \infty,    \\ \nonumber
&&\sum_{n = n_0}^{\pm \infty}|n| \big|a_q^\pm(n) \tilde F^{\pm}(n,n+1) - 
a_q^\pm(n-1) \tilde F^{\pm}(n - 1, n)\big| < \infty.
\eea

(ii) The values $\rho_j \in \R\backslash \Sigma$
are distinct and $\gamma_{\pm, j}\geq 0$ for $1 \leq j \leq q$.

(iii) $\ln(1-|R_\pm|^2)$ is integrable on $\Sigma$ and
$T_\pm(\lam)$ defined via (\ref{rectp}), (\ref{rectm})
extend to single-valued functions on $\Pi_+$ with 
\be \label{hypo 3}
T_-(\lam) = T_+(\lam) \prod_{j=1}^g\frac{\lam - \mu_j^+}{\lam - \mu_j^-}.
\ee

(iv) Transmission and reflection coefficients satisfy (\ref{T E}), (\ref{T infty}), 
and the consistency conditions
\[              
\frac{R_\mp(\lam)}{\ol{R_\pm(\lam)}} = - \frac{T_\pm(\lam)}{\ol{T_\pm(\lam)}}, \qquad
\gamma_{+, j}\, \gamma_{-, j} =  
\frac{\big(\res_{\rho_j} \frac{\prod_{\mu_j^\mp\in M^\mp}(z - \mu_j^\mp)}{\prod_{\mu_j^\pm\in \bar{M}^\pm \cup \tilde{M}^\pm}(z - \mu_j^\pm)} T_\pm(z)\big)^2}
{R_{2g+2}(\rho_j)}.
\]
\end{hypo}

\section{Inverse scattering theory}
\label{secINV}

In this section we reconstruct the operator $H$ from a given set $S_+$ or $S_-$ and 
given quasi-periodic Jacobi operators $H_q^\pm$.

If $S_{\pm}$ (satisfying H.\ref{hypo scat q} (i)--(ii)) and $H_q^\pm$ are known, we
can construct $F^{\pm}(l,m)$ via formula (\ref{glm2 q}) and thus derive the
Gel'fand-Levitan-Marchenko  equation, which has a unique solution by
Theorem~\ref{thmglm}. We obtain that
\begin{align} \nn
K_\pm(n,n) &= \spr{\delta_0}{(1 + \mathcal{F}^{\pm}_n)^{-1} \delta_0}^{1/2}\\
K_\pm(n,n\pm j) &= \frac{1}{K_\pm(n,n)} \spr{\delta_j}{(1 + \mathcal{F}^{\pm}_n)^{-1} \delta_0}.
\end{align}
Since $1 + \mathcal{F}^{\pm}_n$ is
positive, $K_\pm(n,n)$ is positive and we can set (see Theorem~\ref{th:transf})
\begin{align} \label{a+a- q}
a_+(n) &= a_q^+(n)\frac{K_+(n+1, n+1)}{K_+(n,n)}, \\             \nonumber
a_-(n) &= a_q^-(n)\frac{K_-(n, n)}{K_-(n+1,n+1)}, \\             \nonumber
b_+(n) &= b_q^+(n) + a_q^+(n)\frac{K_+(n, n+1)}{K_+(n,n)} - a_q^+(n-1)
\frac{K_+(n-1, n)}{K_+(n-1,n-1)}, \\  \nonumber           
b_-(n) &= b_q^-(n) + a_q^-(n-1)\frac{K_-(n, n-1)}{K_-(n,n)} - a_q^-(n)
\frac{K_-(n+1, n)}{K_-(n+1,n+1)}.
\end{align}
Let $H_+$, $H_-$ be the associated Jacobi operators. As in \cite{emtqps} one proves 

\begin{lemma}             \label{invs1 q}
Suppose a given set $S_{\pm}$ satisfies H.\ref{hypo scat q} (i)--(ii). Then the
sequences defined in (\ref{a+a- q}) satisfy $n |a_{\pm}(n)-a_q^\pm(n)|$,
$n|b_{\pm}(n)-b_q^\pm(n)| \in \ell^1_\pm(\N)$. 

Moreover, $\psi_{\pm}(\lambda, n) = \sum_{m=n}^{\pm \infty}
K_{\pm}(n,m)\psi_q^\pm(\lambda, m)$, where $K_{\pm}(n,m)$ is the 
solution of the Gel'fand-Levitan-Marchenko equation, satisfies 
$\tau_{\pm} \psi_{\pm} = \lambda \psi_{\pm}$.
\end{lemma}

We set 
\be \label{W}
W(\lam)= \frac{R_{2g+2}^{1/2}(\lam)}{T_\pm(\lam)\prod_{j=1}^g (\lam-\mu_j^\pm)}.
\ee
According to H.\ref{hypo scat q}~(iii), this function is holomorphic in
$\C\backslash(M^\pm \cup \bar M^\pm \cup \partial\Sigma )$.

Now we can prove the main result of this section.

\begin{theorem}
Hypothesis H.\ref{hypo scat q} is necessary and sufficient for sets $S_\pm$
to be the left/right scattering data of a unique Jacobi operator $H$ associated with
sequences $a$, $b$ satisfying (\ref{hypo}).
\end{theorem}

\begin{proof}
Necessity has been established in the previous section. By Lemma~\ref{invs1 q}, we
know existence of sequences $a_\pm$, $b_\pm$ and corresponding solutions
$\psi_\pm(z,n)$ associated with $S_+$ (or $S_-$). Hence it remains to establish
$a_+(n)=a_-(n)$ and $b_+(n)=b_-(n)$.

We study the following part of the GLM-equation
\be
\sum_{m \in \Z} \Phi_+(n,m)\ol{\psi_q^+(\lam,m)}, \qquad 
\Phi_+(n, .) := \sum_{l=n}^{\infty}K_+(n,l)
{\tilde F}^+(l,.) \in \ell_+^1(\Z)
\ee
and obtain as in \cite[Theorem 8.2]{emtqps}
\be      \label{h3 q}
T_+(\lam)h_-(\lam, n) = \ol{\psi_+(\lam, n)} 
+ R_+(\lam) \psi_+(\lam, n), \qquad \lam \in \Sigma,
\ee
where
\begin{align} \nonumber         \label{h_mp(w, n)}
h_-(\lam, n) &= \frac{\ol{\psi_q^+(\lam,n)}}{T_+(\lam)}\bigg( 
\frac{1}{K_+(n, n)} + 
\sum_{m = - \infty}^{n-1}\Phi_+(n, m) 
\frac{\ol{\psi_q^+(\lam,m)}}{\ol{\psi_q^+(\lam,n)}}  \\
& \quad + \sum_{j=1}^q \gamma_{+, j} \hat \psi_+(\rho_j, n)  
\frac{W_{n-1}(\hat \psi_q^+(\rho_j), \ol{\psi_q^+(\lam)})}
{(\lam - \rho_j) \ol{\psi_q^+(\lam,n)}}\bigg).             
\end{align}
Clearly, using $\bpsi_q^\pm(\lam,n)=\ol{\psi_q^\pm(\lam,n)}$, the function $h_+(\lam,n)$
extends to a meromorphic function in $\C\backslash\Sigma$.

Similar we obtain a function $h_+(z,n)$. The functions $h_{\mp}(z, n)$ have  
simple poles at $\mu_j^\mp \in M^\mp$ and the following limits
\begin{align} \label{h infty}
h_{\mp}(z, n) &= \frac{z^{\pm n}}{K_\pm(n,n)T_\pm(\infty)} 
\Big(\prod_{j=0}^{n-1}a_q^\pm(j)\Big)^{\mp1}\Big(1+O(\frac{1}{z})\Big), 
\quad z \rightarrow \infty, \\ \label{h(rho)}
\lim_{z\rightarrow \rho_j}h_{\mp}(z, n) &= 
\frac{\gamma_{\pm,j}\hat \psi_\pm(\rho_j,n)}
{\sqrt{\gamma_{+,j}\gamma_{-,j}}\prod_{\mu \in M^\mp}(\rho_j - \mu)}.
\end{align}
By virtue of the consistency condition $T_\pm(\lam) \ol{R_\pm(\lam)} = - \ol{T_\pm(\lam)} R_\mp(\lam)$ 
we obtain
\begin{equation} \label{h4 q}
\overline{h_{\pm}(\lam, n)} + R_{\pm}(\lam)h_{\pm}(\lam, n) = 
\psi_{\mp}(\lam, n) T_\pm(\lam), \qquad \lam \in \Sigma.
\end{equation}
Eliminating $R_{\pm}(\lam)$ in (\ref{h3 q}) and (\ref{h4 q}) and using (\ref{W}) yields
\begin{align} \nn
&R_{2g+2}^{-1/2}(\lam) \prod_{j=1}^g (\lam-\mu_j^\pm)
\big(\ol{h_\pm(\lam,n)} \psi_\pm(\lam,n) 
- \ol{\psi_\pm(\lam,n)} h_\pm(\lam,n)\big)   
\\       \label{def G}  
&=\frac{\psi_+(\lam,n) \psi_-(\lam,n) -h_+(\lam,n)h_-(\lam,n)}{W(\lam)}
= G(\lam,n), \qquad \lam \in \Sigma. 
\end{align}
Observe that $G(\lam, n)$ extends to a meromorphic function
on $\C\backslash\Sigma$ which is continuous on the interior
of $\Sigma$ with equal real limits from above and below:
$G(\lam\pm\I0, n) = \overline{G(\lam\mp\I0, n)} = G(\lam\pm\I0, n)$.
So by the Schwarz reflection principle $G(z, n)$ is holomorphic on the set 
$\C\backslash(M^\pm \cup \bar M^\pm \cup \{\rho_k\}_{k=1}^q \cup \partial\Sigma)$. 

Since the difference $\psi_+ \psi_- - h_+ h_-$ vanishes at the 
points $\rho_k$ by (\ref{h(rho)}) and H.\ref{hypo scat q}~(iv),
the poles there are removable. 

Next let us investigate the behavior at the band edges. 
Since $G(z, n)$ has at most poles, it is sufficient to control the behavior
from one direction. To this end we use (\ref{T E}) and the identity (\ref{def G}). If 
$E \notin \ti M^+ \cup \ti M^-$, 
\begin{align} \nonumber
&\lim_{\lam \rightarrow E}  R_{2g+2}^{1/2}(\lam)  \prod_{j=1}^g(\lam-\mu_j^\pm)
h_\mp(\lam, n) \overline{\psi_\mp(\lam,n)}\\ \nonumber
&=\lim_{\lam \rightarrow E} \frac{R_{2g+2}^{1/2} \prod_{j=1}^g(\lam-\mu_j^\pm)}{T_\pm}
\Big( \overline{\psi_\pm} + R_\pm \psi_\pm \Big) \overline{\psi_\mp}
\\                              \nonumber
&=\lim_{\lam \rightarrow E} \frac{R_{2g+2}^{1/2} \prod_{j=1}^g(\lam-\mu_j^\pm)}{T_\pm}
\Big( (R_\pm + 1) \psi_\pm + \overline{\psi_\pm} - \psi_\pm \Big) 
\overline{\psi_\mp} = 0. 
\end{align} 
For $E \in \ti M^+ \cap \ti M^-$, the functions $\psi_\pm(\lam, .) + \ol{\psi_\pm(\lam, .)}$
are bounded in the vicinity of $E \in \ti M^\pm$,
\begin{align*}  
&\lim_{\lam \rightarrow E} \frac{R_{2g+2}^{1/2} \prod_{j=1}^g(\lam-\mu_j^\pm)}{T_\pm}
\Big( \overline{\psi_\pm} + R_\pm \psi_\pm \Big) \overline{\psi_\mp} \\  
&=\frac{R_{2g+2}^{1/2} \prod_{j=1}^g(\lam-\mu_j^\pm)}{T_\pm}
(\psi_\pm + \ol{\psi_\pm}) \overline{\psi_\mp} + \prod_{j=1}^g(\lam-\mu_j^\pm)\psi_\pm \ol{\psi_\pm}
\frac{R_{2g+2}^{1/2} (R_\pm - 1)}{T_\pm} = 0.
\end{align*} 
Analogously, one proves that $R_{2g+2}(z)G(z,n)$ vanishes at the points 
$E \in \ti M_\pm \backslash \ti M_\mp$. Therefore this function is
holomorphic on $\C$ and $R_{2g+2}(z)G(z,n)=O(z - E)$ for $E \in \partial \Sigma$.  
Thus $G(z,n)$ has only removable singularities and since $G(z,n) \rightarrow 0$
as $z \to \infty$, Liouville's theorem implies $G(z,n) \equiv 0$ and
\[
\psi_+(z,n) \psi_-(z,n)- h_+(z,n)h_-(z,n) \equiv 0, \qquad \forall n \in \Z.
\]
For $z \rightarrow \infty$ we obtain by (\ref{B4jost}) and (\ref{h infty}) 
\[
T_+(\infty)T_-(\infty)=\bigg(\frac{1}{K_+(0,0)K_-(0,0)}\prod_{j=0}^{n-1}\frac{a_-(j)}{a_+(j)}\bigg)^2
\]
and therefore by (\ref{T infty})
\be
a_+(n) = a_-(n) \equiv a(n), \qquad \forall n \in \Z.
\ee
It remains to prove $b_+(n)=b_-(n)$. Proceeding as for $G(\lam,n)$ we can show
that
\begin{align}    \nonumber
&\frac{\psi_+(\lam,n) \psi_-(\lam,n+1) 
-h_+(\lam,n+1) h_-(\lam,n)}{W(\lam)}\\
&= \frac{\prod_{j=1}^g (\lam-\mu_j^+)}{R_{2g+2}^{1/2}(\lam)}
\big(\overline{h_+(\lam,n+1)} \psi_+(\lam,n) 
- \overline{\psi_+(\lam,n)} h_+(\lam,n+1)\big)
\end{align}
is a constant equal to $-1/a(n)$. Thus
\[
\bar W(z,n) := a(n) \left(\psi_+(z,n) \psi_-(z,n+1) - h_+(z,n+1) h_-(z,n) \right)\\
= - W(z).
\]
Computing the asymptotics as $z \rightarrow \infty$ (compare (\ref{B4jost})) we see
\be
0= \bar W(z,n) - \bar W(z,n-1) = \frac{b_+(n) - b_-(n)}{A_+(0)A_-(0)}
\ee
and in particular $b_+(n) = b_-(n) \equiv b(n)$.
\end{proof}

{\bf Acknowledgments.}
I.E. gratefully acknowledges the extraordinary hospitality of the
Faculty of Mathematics of the University of Vienna during two stays in 2006,
where parts of this paper were written.

\end{document}